\documentstyle[11pt]{article}
\title{One page proof of the Riemann hypothesis}
\author{Andrzej M\c{a}drecki\thanks{Institute of Mathematics and
Computer Science, Wroc{\l}aw University of Technology(WUT),
51-200 Wroc{\l}aw, Poland}}
\newtheorem{th}{Theorem}
\newtheorem{pr}{Proposition}
\newtheorem{re}{Remark}
\newfont{\lll}{msbm10 scaled 1095}

\def\LC{\mbox{\lll \char67}}
\def\LM{\mbox{\lll \char77}}
\def\LN{\mbox{\lll \char78}}
\def\LP{\mbox{\lll \char80}}
\def\LR{\mbox{\lll \char82}}
\def\LZ{\mbox{\lll \char90}}
\begin{document}
\maketitle

{\bf Abstract}. We give a short Wiener measure proof of the Riemann
hypothesis based on a surprising, unexpected and deep relation between
the Riemann zeta $\zeta(s)$ and the trivial zeta $\zeta_{t}(s):=Im(s)(2Re(s)-1)$.

\section{Functional analysis and probability theory of the Riemann hypothesis.}

Let $C^{+} = C^{+}(\LR)$ be the real vector space of all real valued, {\bf symmetric}
and {\bf continuous} functions defined on the {\bf real numbers field} $\LR$.
We also consider the following {\bf normed spaces}, {\bf Banach spaces} and
{\bf Frechet spaces} . 

By $S_{2}$ we denote the real vector subspace of functions $c \in C^{+}$
with the finite {\bf second-max-moment}, i.e.
\begin{displaymath}
\mid\mid c \mid\mid_{2}:=maximum_{x \in \LR}\mid x^{2} c(x) \mid <+\infty,
\end{displaymath}
whereas by $S_{0}$ we denote the real vector subspace of $C^{+}$ of all
{\bf bounded continuous functions} with the sup-norm $\mid\mid \cdot
\mid\mid_{0}$. Obviously, $B_{2}:=(S_{2}, \mid\mid \cdot \mid\mid_{2})$ is
a normed space and $B_{0}:=(S_{0}, \mid\mid \cdot \mid\mid_{0})$ is a Banach
space.

Let us denote by $I$ the unit interval $[0,1]$ and by $I^{c}=(1,+\infty)$
its complement. Moreover, $L^{1}(I^{c}, dx)$ is the Lebesgue space of
all real absolute integrable functions $f : I^{c}\longrightarrow \LC$
with the {\bf finite first integral moment}
\begin{displaymath}
\mid\mid f \mid\mid_{1}:=\int_{I^{c}}\mid f(x) \mid dx<+\infty,
\end{displaymath}
and $dx$ is the {\bf Lebesgue measure}. In particular, since
\begin{displaymath}
\int_{I^{c}}\mid f(x) \mid dx \le \mid\mid f
\mid\mid_{2}\int_{1}^{\infty}\frac{dx}{x^{2}},
\end{displaymath}
we have  $S_{2}\mid I^{c} \subset L^{1}(I^{c},dx)$.

We consider the {\bf canonical Fourier (cosine) transform} ${\cal
F}:S_{2}\longrightarrow S_{0}\cap C^{+}$ defined by the well-known formula :
\begin{equation}
({\cal F}f)(x)\;:=\;\int_{\LR}e^{2\pi
ixy}f(y)dy=2\int_{0}^{+\infty}cos(2\pi xy)f(y)dy =:\hat{f}(x),
\end{equation}
where $x \in \LR_{+}:=[0,+\infty)$ and $f \in S_{2}$.

The formula
\begin{equation}
\theta(f)(x):=\sum_{n=1}^{\infty}f(nx)=\int_{\LN^{*}}f(nx)dc_{\LZ}(n);x>
0, f \in S_{2},
\end{equation}
defines the canonical {\bf Jacobi theta transform} $\theta : S_{2}
\longrightarrow C^{+}(\LR-\{0\})$, where $\LN^{*},\LN$ and $\LZ$ are the
{\bf multiplicative semigroup of positive integers}, {\bf additive semigroup
of positive integers} and the {\bf ring of integers}, respectively and
$c_{\LZ}$ marks the {\bf calculating (Haar) measure} of $(\LZ, +)$ normalized
as : $c_{\LZ}(\{0\})=1$.

Finally, by $M : S_{2}\longrightarrow C(0,3]\otimes_{\LR}\LC $ we denote the
{\bf Mellin transform}, i.e.
\begin{equation}
M(f)(s)\;:=\;\int_{0}^{\infty}x^{s}f(x)\frac{dx}{x}\;Re(s)>0, f\in S.
\end{equation}

All in the sequel $\LC$ marks the {\bf complex number field} and if $s
\in \LC$ is arbitrary then by $Re(s)$ and $Im(s)$ we mark the {\bf real}
and {\bf imaginary} part of $s$.

Obviously, the transforms ${\cal F}, \theta$ and $M$ are {\bf continuous
operators} defined on the Banach space $B_{2}$ with values in suitable
Frechet spaces, according to the inequalities:
\begin{displaymath}
\sum_{n=1}^{\infty}\mid f(nx)\mid \le \frac{\mid\mid f \mid\mid_{2}\zeta(2)}
{x^{2}},
\end{displaymath}
and
\begin{displaymath}
\mid M(f)(s)\mid \le \mid\mid f \mid\mid_{0}\int_{0}^{1}\mid x^{s-1}\mid dx
+\mid\mid x^{s-2}f(x) \mid\mid_{0}\int_{1}^{\infty}\frac{dx}{x^{2}}.
\end{displaymath}

In the sequel we work with the fundamental {\bf Poisson cylinder}
\begin{equation}
\LP:=\{c \in C^{+}:\;c(0)=1\;and\; \exists(f \in S_{2})with(\hat{f}\in
S_{2})(c\;=\;f+\hat{f})\}
\subset S_{2}+S_{2}=S_{2}.
\end{equation}
In particular, for each $c \in \LP,  {\cal F}(c)=c$ i.e. $\LP$ is a subcylinder
of the eigenvectors space of ${\cal F}$, which corresponds to the
eigenvalue $+1$ , i.e. $\LP$ is the set of {\bf fixed points of} ${\cal F}$.

It is well-known that the $\LR$-vector space $C^{+}$ has got the natural
structure of the {\bf Frechet space}.

Let ${\cal P}$ be the $\sigma$-{\bf field} of all subsets of $\LP$
generated by {\bf cylinders} $C$ of $\LP$ of the form:
\begin{displaymath}
C\;=\;C(t_{1}, ... , t_{n};B):=\{c \in \LP : (c(t_{1}), ... ,
c(t_{n}))\in  B\},
\end{displaymath}
where $B$ is a {\bf Borel subset} from $\LR^{n}$ and $t_{0}=0<t_{1}< ...
<t_{n}$.

The importance of the {\bf phase space} $(\LP, {\cal P})$ is motivated
by the two facts :

(1) the functions from the Poisson cylinder $\LP$ satisfies the
fundamental in the zeta theories {\bf Poisson Summation Formula}(PSF in
short), and

(2) On ${\cal P}$, there exists fundamental for this short Wiener measure proof
- the {\bf non-trivial Wiener-Riemann measure} $r : {\cal P}\longrightarrow
[0,1]=:I$.

Our method of transferring results concerning the functional analysis
and probability theory of RH looks, broadly speaking, as follows :
assume that we have the classical Riemann continuation equation (1.34)
and the Riemann hypothesis.

Now having the analytic number theory problem of RH, we can try to
solve the corresponding functional analysis and probabilistic problem
for the extension of the above mentioned functional equation, for the
functional space $\LP$, which may be easier than in the original RH
setting, since we have additional functional analysis and probability
means at our disposal. Having done this, we can again try to put those
functional-probability solution together in some way and this may
happen to yield a solution of the RH-problem.

A {\bf stochastic process} $B = (B_{t} : t \ge 0)$ defined on a
probability space $(\Omega, {\cal A}, Prob)$ is said to be the {\bf
standard Brownian motion} iff it is {\bf gaussian} (i.e. its all
1-dimensional {\bf distributions are gaussian}), its {\bf moment
function} is zero : $EB_{t} = 0$, and its {correlation function}
$EB_{t}B_{s}$ is equal to $min(t,s)$.

Here and all in the sequel $EX$ marks the {\bf expected value} of a
{\bf real random variable} $X$ (rv for short).

For the existence of the Brownian motion see e.g. [W, II.3]. In particular
, $B$ has {\bf continuos paths} (since it satisfies the well-known {\bf
Kolmogorov condition} - see [W, II.4 , Th.4.5]) and the {\bf
independent increments}. Finally $B$ gives the main example of a {\bf Markov
process} and a {\bf martingale}.

Let us now consider the {\bf peak function} $p(t)$ defined like :
\begin{equation}
p(t)\;=\;1\;-\; t \;\;if\;\;0 \le t  \le 1,
\end{equation}
and
\begin{displaymath}
p(t)\;=\;0\; if\;\; t  \ge 1,
\end{displaymath}
(in particular, $p$ is in Cameron-Martin space) and the stochasic process
\begin{equation}
B_{t}^{p}(\omega):=\;(B_{t}(\omega)\;+\;p(t))\;\;;t \ge 0, \omega \in
\Omega,
\end{equation}
as a {\bf random element} (re in short) $B^{p} : (\Omega, {\cal A},
Prob)\longrightarrow (C^{+}, {\cal C}^{+})$ , where ${\cal C}^{+}$ is
the {\bf cylinder} $\sigma$- {\bf field} of $C^{+}$.

Since $B$ determines the {\bf standard Wiener measure} $w$ on the phase
space $(C^{+}, {\cal C}^{+})$ by the well-known formula (it is the {\bf law}
of $B$) :
\begin{equation}
w(C)\;:=\;Prob(B^{-1}(C)), \;\; C \in {\cal C}^{+},
\end{equation}
then on the Poisson phase space $(\LP, {\cal P})$ we can define the following
{\bf Wiener-Riemann measure} $r$ according to the formula :
\begin{equation}
r(P)\;\;:=\;\;\sum_{n=1}^{\infty}\frac{1}{2^{n}}Prob(\omega \in \Omega:
G(t)B_{\sqrt{t}}^{p}(\omega)(t)\in P, t \in [n-1,n] )\;=\;
\end{equation}
\begin{displaymath}
\sum_{n=1}^{\infty}\frac{1}{2^{n}}w((G^{-1}P-p)\cap
e_{n}(C[\sqrt{n-1},\sqrt{n}])),\;P \in {\cal P},
\end{displaymath}
where

(1) $G(t) = e^{-\pi t^{2}}$ is the {\bf standard Gauss function},

(2) $C[\sqrt{n-1},\sqrt{n}]$ is the Banach space of all real valued continuous functions
defined on the segment $[\sqrt{n-1},\sqrt{n}]$ and considered as the subspace of $C^{+}$
through the {\bf embedding} $e_{n} : C[\sqrt{n-1},\sqrt{n}]\longrightarrow C^{+}$ by the
formula : $e_{n}(c)(t) :=c(t)$ if $\sqrt{n-1} \le t \le \sqrt{n}$;
$e_{n}(c)(t):=c(\sqrt{n-1})$ if $t \le \sqrt{n-1}$ and $e_{n}(c)(t):=c(n)$ if
$t \ge \sqrt{n}$.
\begin{re}
Let us observe that $e_{n}$ is a linear monomorphism , i.e.
$Ker(e_{n})=\{0\}$. Thus although $C[\sqrt{n-1},\sqrt{n}]$ is not {\bf formally} a
subspace of $C^{^{+}}$ (since the product $C[\sqrt{n-1},\sqrt{n}] \cap C^{+}$ is
empty), however we can identify here $C[\sqrt{n-1},\sqrt{n}]$ with its isomorphic image
$Im(e_{n})=e_{n}(C[\sqrt{n-1},\sqrt{n}])$.

Let us also observe, that exactly therefore since the space (cone) $\LP$
of functions which are self-similar with respect to Fourier transform
is a priori "small" for Wiener measure then a posteriori the Wiener-Riemann
measure $r$ is the {\bf series} of measures $\{r_{n}\}$ being the restriction to $\{C[\sqrt{n-1},
\sqrt{n}]\}$ of the measure $(G^{-1})^{*}w_{p}$, which is subsequently
the transport by $G^{-1}$ of the shifted Wiener measure by $p$ : $w_{p}(X):=
w(X+p)$ (which is {\bf equivalent} to $w$, since $p$ belong to the Cameron-
Martin space).
\end{re}

\begin{pr}({\bf On the existence of the Wiener-Riemann measure and its
RH-properties}.)

(I). For each $c \in \LP$  and $x>0$ the {\bf Poisson Summation Formula}(PSF
in short) holds, i.e.
\begin{equation}
\frac{1}{x}\theta(c)(\frac{1}{x})\;+\;(c(0)=1/2x) \;=\;1/2\;+\;
\theta(c)(x).
\end{equation}
(II). The measure space $(\LP, {\cal P}, r)$ has the following four
properties :

($r_{0}$)({\bf Non-triviality}). $r(\LP) = 1$.

($r_{1}$)({\bf Starting point}). $\int_{\LP}c(0)dr(c) = p(0) =1$.

($r_{2}$)({\bf Vanishing of moments}). For all $t \ge 1$ holds
\begin{displaymath}
\int_{\LP}c(t)dr(c)\;=\;0.
\end{displaymath}
($r_{3}$)({\bf The Fubini obstacle-Hardy-Littlewood theorem obstacle -
Existence of moments of 1/2-stable Levy distributions}).

The double integrals (the averaging of the Mellin transform w.r.t. $r$)
\begin{equation}
b_{s}\;:=\;\int \int_{\LP \times \LR_{+}}\mid x^{s-1}c(x) \mid
dr(c)dx=\int_{\LP}M(\mid c \mid)(re(s))dr(c) \;=\;
\end{equation}
\begin{displaymath}
\;=\;\sum_{n=1}^{\infty}(1/2^{n})E[M(\chi_{[\sqrt{n-1},\sqrt{n}]}G^{-1}B^{p}
_{\sqrt{\cdot}})](re(s)),
\end{displaymath}
are

(i) {\bf finite} if $Re(s) \in (0,1/2)$, and 

(ii){\bf infinite} if $Re(s) \ge 1/2$.
\end{pr}

{\bf Proof}. (I). It is widely known fact among specialists on zetas. We
remark at once that the form (1.9) of (PSF) is a consequence of the fact that each
$c \in \LP$ is a {\bf fixed point} of ${\cal F}$, i.e. ${\cal F}(c) = c$.
Let $c \in \LP$ be arbitrary. Then $c$ is a continuous function in
$L_{1}(\LR)$. We show that $c$ satisfies the following two conditions :

(i) The series $\sum_{n \in \LN}c(x+n)$ is {\bf uniformly convergent} for
all $x \in (0,1)=:D$.

Reely, it is nothing that $\theta(c_{+x})(1)$ and obviously
\begin{displaymath}
max_{x \in D} \mid \sum_{n \ge N}c(x+n) \mid \le max_{x \in \LR}\mid
(x+n)^{2}c(x+n) \mid \sum_{n \ge N}max_{x \in D}(\frac{1}{(x+n)^{2}})
\le  \mid\mid c \mid\mid_{2}\zeta(2),
\end{displaymath}
so, according to the Weierstrass critirion, the above series are
uniformly convergent in $D$. In particular 

(ii) $\theta(c)(1)=\sum_{n \in \LN^{*}}c(n)$ is {\bf convergent}.

One of the wider formulations of (PSF) and its proof, in the general
case of any LCA groups, a reader can find in the beautiful {\bf
Narkiewicz's book} [N, Appendix I. 5, Th. VIII]. We used above Th. VIII
in the case : $G= \LR, H=K=\LZ$ and $D=(0,1)$.

Now let $x>0$ be arbitrary. We apply the above (PSF) in the case of
the function
\begin{displaymath}
c_{\cdot x}(y) \;:=\;c(xy)\;,\;c \in \LP.
\end{displaymath}
Since obviously, we have the following easy  calculus :
\begin{displaymath}
\hat{c}_{\cdot x}(y)=2\int_{0}^{\infty}cos(2\pi
yz)c(xy)dz=\{u:=xz\}=(2/x)\int_{0}^{\infty}cos(2\pi \frac{yu}{x})c(u)du
=
\end{displaymath}
\begin{displaymath}
=\frac{1}{x}\hat{c}{y/x}\;=\;\frac{1}{x}c(\frac{y}{x}).
\end{displaymath}
Hence
\begin{displaymath}
\theta(c)(x)\;+\;\frac{1}{2}\;=\;\frac{1}{2x}\;+\;\frac{1}{x}\theta(c)(
\frac{1}{x}).
\end{displaymath}
$(II)(r_{0})$. Let us observe that for each $n \in \LN^{*}$  :
\begin{displaymath}
\LP \cap C[\sqrt{n-1}, \sqrt{n}]\;=\;C[\sqrt{n-1},\sqrt{n}].
\end{displaymath}
Really, let $f \in C[\sqrt{n-1},\sqrt{n}]$ be arbitrary. We can consider $f$ like a
restriction of a function $\tilde{f}$ from $C^{+}$ with a {\bf
compact support}. Let us consider the {\bf second order Fredholm integral
equation} of the form (Fox equation) :
\begin{equation}
f(x)\;=\;g(x)\;+\;2\int_{0}^{\infty}cos(2\pi xy)g(y)dy .
\end{equation}
In [KKM, II. 23, Example] the following Fredholm-Fourier-Fox was considered
integral equation, with a {\bf parameter} $\lambda$ :
\begin{equation}
f(x)\;=\;\phi(x)\;-\;\lambda
\sqrt{\frac{2}{\pi}}\int_{0}^{\infty}\phi(x)cosxt dt.
\end{equation}
The authors solved that Fredholm equation by using the {\bf Mellin transform}
and they have obtained the following formula for the solution :
\begin{equation}
\phi(x)\;=\;\frac{f(x)}{1-\lambda^{2}}\;+\;\frac{\lambda}{1-\lambda^{2}}
\sqrt{\frac{2}{\pi}}\int_{0}^{\infty}f(t)cosxt dt .
\end{equation}
Thus taking $\lambda = -\sqrt{\pi/2} \ne \pm 1$ in (1.12) and (1.13) , we
obtain the existence of $g$ with the properties :
\begin{displaymath}
g(x)\;=\;\frac{\pi f(x)}{\pi - 2}\;+\;\frac{\pi \hat{f}(x)}{\pi - 2},
\end{displaymath}
and such that for $x \in [\frac{\sqrt{(n-1)}}{2\pi}, \frac{\sqrt{n}}{2\pi}]$
holds
\begin{displaymath}
f(2\pi x)\;=\;g(2\pi x)\;+\;\hat{g}(2\pi x).
\end{displaymath}
Let us observe that we can put $f(0)=1$ since $supp(w_{p})=C_{1}[0,1]$
(or equivalently $B^{p}_{0}=1$ with probability 1). Moreover $g \in S_{2}$
since we have the following trivial calculus :

\begin{displaymath}
max_{x \in \LR}\mid x^{2}F(x) \mid \le max_{[0,1]}\mid
x^{2}F(x)\mid\;+\; max_{x \ge 1}\mid x^{2}F(x)  \mid,
\end{displaymath}
for any continuous function $F(x)$. Making the substitution : $u=x^{2}t$
in the integral $x^{2}\int_{0}^{n(F)}F(t)coxt dt$ where $n(F)>0$ is
such that $supp(f) \subset [0,n(F)]$ we obtain the following easy
estimation :
\begin{displaymath}
\mid x^{2}\int_{0}^{n(F)}F(t)cosxt dt \mid \le max_{u \ge 0}\mid
F(\frac{u}{x^{2}})\mid,
\end{displaymath}
which is finite and do not depends on $x$ , for $x\ge 1$.

Hence
\begin{displaymath}
r(\LP)=\sum_{n=1}^{\infty}\frac{1}{2^{n}}Prob(G(\cdot)B^{p}_{\sqrt(\cdot)}
\in \LP \cap e_{n}C[n-1,n])=\sum_{n=1}^{\infty}\frac{w(e_{n}C[n-1,n])}{2^{n}}
=1.
\end{displaymath}
($r_{1}$). We have
\begin{displaymath}
\int_{\LP}c(0)dr(c)=E(B_{0}+p(0))=p(0)=1.
\end{displaymath}
($r_{2}$). If $t \ge 1$ then $t \in [n-1,n]$ for some $n \ge 2$ and
\begin{displaymath}
\int_{\LP}c(t)dr(c)=G(t)EB_{\sqrt{t}}/2^{n}=0.
\end{displaymath}
$(r_{3})$. This is the unique non-trivial property of $r$. The idea of
this proof is taken from the Proposition 3 of [AM].

We first give the {\bf upper estimation} of the iterated integral
\begin{displaymath}
I_{dxdr}(s)\;:=\;\int_{0}^{\infty}dx(\int_{\LP}\mid x^{s-1}c(x)\mid
dr(c)),
\end{displaymath}
for $s$ with $u=Re(s) \in (0,1/2)$.

Let us observe that
\begin{equation}
\int_{\LP}\mid c(x) \mid dr(c)\;=\;\frac{G(x)E\mid
B^{p}_{\sqrt{x}}\mid}{2^{n}},
\end{equation}
if $n-1 \le x \le n, n=1,2, ...$. Hence
\begin{equation}
\int_{0}^{\infty}x^{u-1}dx \int_{\LP}\mid c(x) \mid dr(c)\le
\int_{0}^{1}x^{u-1}G(x)p(x)dx+\sum_{n=2}^{\infty}\frac{1}{2^{n}}\int_{n-1}
^{n}x^{u-1}G(x)E\mid B_{\sqrt{x}} \mid dx\le
\end{equation}
\begin{displaymath}
\le \int_{0}^{1}x^{u-1}G(x)p(x)dx+\int_{0}^{\infty}x^{u-1}G(x)E\mid
B_{\sqrt{x}}\mid dx.
\end{displaymath}

Since in the inequality (1.15) the first unproper or proper integral obviously
exists, thus the problem of the {\bf convergence} of $I_{dxdr}(s)$ is reduced to the
convergence of the integral
\begin{equation}
b_{s}:=\int_{0}^{\infty}x^{u-1}G(x)E\mid B_{\sqrt{x}} \mid
dx=1/\sqrt{2\pi}\int_{0}^{\infty}x^{u-1}G(x)(\int_{\LR}\mid y+y_{0} \mid
e^{-\frac{(y+y_{0})^{2}}{2x}}/\sqrt{x}dy)dx,
\end{equation}
for any $y_{0}>0$, according to the facts that the distribution of
$B_{\sqrt{x}}$ is gaussian with {\bf mean zero} and {\bf variance}
$\sqrt{x}$ and the {\bf translational invariance} of the {\bf Lebesgue measure}
$dy$.

Now, let us observe (what was first observed in [AM]), that the iterated
integral in the right-hand side of (1.16) - according to :
\begin{equation}
e^{-y^{2}/2x}e^{-2\mid y \mid y_{0}/2x}\;\le\;e^{-\mid y \mid y_{0}/x}
\end{equation}
and a suitable substitution can be approximated as follows :
\begin{equation}
\le (\int_{0}^{\infty}x^{u}\frac{e^{-y_{0}^{2}/2x}}{\sqrt{2
\pi}x^{3/2}}dx)(y_{0}^{-2}max_{x \ge 0}(x^{2}G(x))\int_{\LR}\mid t \mid
e^{-\mid t \mid}dt+max_{x \ge 0}(xG(x))\int_{\LR}e^{-\mid t \mid}dt) .
\end{equation}

But now, the function
\begin{equation}
d_{y_{0}}(x)\;:=\;\frac{y_{0}e^{-y_{0}^{2}/2x}}{\sqrt{2\pi}x^{3/2}}\;,\;
x>0,
\end{equation}
is exactly the {\bf density} of a rv $L_{y_{0}}$ with the {\bf
1/2-stable-Levy distribution} (with a parameter $y_{0}$). It is well-known
that the distribution of $L_{y_{0}}$ is concentrated on $\LR_{+}$ and
that it is the {\bf unique p-stable distribution} with $p \in (0,1)$, which
has an elementary analytic and simple formula for the density of a power
-exponential form (see e.g. [AM]).

The most important fact concerning $L_{y_{0}}$, what we use for this short
proof of (RH) is the problem of the {\bf existence of the Orlicz moments}
of $L_{y_{0}}$ (there are also Hilbert moments $EX^{2}$ and Banach
moments $E\mid X \mid^{p}, p\ge 1$).
More exactly, it is well-known that (see [F]):
\begin{equation}
E(L_{y_{0}}^{u})
\;=\;\int_{0}^{\infty}x^{u}d_{y_{0}}(x)dx<+\infty\;if\;u \in (0,1/2),
\end{equation}
and
\begin{equation}
E(L_{y_{0}}^{u})\;=\;+\infty\;if\;u \ge 1/2.
\end{equation}

Combining (1.18) with (1.20), we finally obtain that the iterated integral $I_{dxdr}(s)$
is {\bf finite} if $Re(s) \in (0,1/2)$. Since obviously the measures $r$ and
$dx$ are $\sigma$-{\bf finite}, then according to the {\bf Tonelli-Fubini theorem}
(TF in short), the numbers $b_{s}$ are {\bf finite} in this case.

The below lower estimation - also modeled on the previous one - shows that
(TF) is {\bf violated} in the case of the triplet :
\begin{displaymath}
(b_{s}=\int \int_{\LP \times \LR}\mid x^{s-1}c(x)\mid
dr(c)dx,I_{dxdr}(s), I_{drdx}(s))
\end{displaymath}
for $Re(s) \ge 1/2$. Thus the violation of (FT) in the above case is mainly
{\bf responsible} for the {\bf non-triviality} of the {\bf Riemann hypothesis}
( the {\bf Hardy-Littlewood theorem}: on the critical line $Re(s)=1/2$
the Riemann zeta $\zeta(s)$ has {\bf infinitely many zeros}).

We have
\begin{equation}
b_{1/2}:= \int_{0}^{\infty}\frac{dx}{\sqrt{x}}(\int_{\LP}\mid c(x) \mid
dr(c))\ge
\sum_{n=1}^{\infty}\frac{1}{2^{n}}\int_{n-1}^{n}\frac{G(x)E\mid
B^{p}_{\sqrt{x}}\mid dx}{\sqrt{x}}\ge
\end{equation}
\begin{displaymath}
\ge
\sum_{n=2}^{\infty}\frac{G(x_{n})}{2^{n}}\int_{n-1}^{n}\frac{dx}{\sqrt{x}}
\int_{\LR}\frac{\mid y \mid e^{-y^{2}/2x}dy}{\sqrt{2\pi
x}}=:\sum_{n=2}^{\infty}g_{n}I_{n} \ge
\end{displaymath}
But using the classical {\bf Tshebyshev inequality} for positive monotonic
finite real sequences (see [Mi, I. 9]) - for each $N \ge 2$ we obtain
\begin{equation}
(\sum_{n=2}^{N+1}g_{n}I_{n})\;\ge(\frac{1}{N}\sum_{n=2}^{N+1}g_{n})(\sum
_{n=2}^{N+1}I_{n}),
\end{equation}
where
\begin{equation}
(\sum_{n=2}^{N+1}I_{n})\;=\;(2\pi)^{-1/2}\int_{\LR}dy \mid y
\mid(\int_{1}^{N+1}\frac{e^{-y^{2}/2x}dx}{x}).
\end{equation}
Combining (1.22), (1.23) and (1.24) we claim that for all $p \ge 1$ and $N \ge
2$ holds :
\begin{equation}
b_{1/2}^{1/p} \ge (2\pi)^{-1/2p}(\frac{\sum_{n=2}^{N+1}g_{n}}{N})^{1/p}
(\int_{\LR}dy \mid y \mid
\int_{1}^{N+1}\frac{e^{-y^{2}/2x}dx}{x})^{1/p}.
\end{equation}
But it is well-known fact (see e.g. [Mu, II . 5, Exercises 1 and 2])
that from the {\bf Holder inequality} follows that for any $g=(g_{2}, ...
,g_{N+1})\ge (0, ... ,0)$ the function
\begin{displaymath}
f(p)\;=\;(\frac{1}{N}\sum_{n=2}^{N+1}g_{n})^{1/p}
\end{displaymath}
is {\bf non-decreasing} and {\bf bounded} in $p$ and moreover
\begin{equation}
\lim_{p \rightarrow
\infty}(\frac{1}{N}\sum_{n=2}^{N+1}g_{n})^{1/p}=max_{2\le n \le
N+1}(g_{n}).
\end{equation}
Combining (1.25) with (1.26) and observing that for each $y \in \LR$ the
minimum below is not zero:
\begin{equation}
min_{1 \le x \le N+1}e^{-y^{2}/2x}\;=\;e^{-y^{2}/2}>0 ,
\end{equation}
we obtain
\begin{equation}
\lim_{p}b_{1/2}^{1/p}\ge max_{2 \le n \le
N+1}(g_{n})(=G(x_{2})/4)\lim_{p}(\int_{\LR}\mid y \mid
e^{-y^{2}/2}dy)^{1/p}(\int_{1}^{N+1}\frac{dx}{x})^{1/p},
\end{equation}
for each $N \ge 2$. Since the left hand side {\bf does not depend on}
$N$, then we finally get
\begin{equation}
\lim_{p}b^{1/p}_{1/2}\;\ge
\;\frac{G(x_{2})}{4}\lim_{p}(\int_{1}^{\infty}\frac{dx}{x})^{1/p}=+\infty.
\end{equation}
So, it must be that $b_{1/2}= \infty$ and the Fubini-Tonelli theorem is
violated for $Re(s)=1/2$.

\begin{pr}({\bf The Muntz relations for $(\zeta(s), [s(s-1)]^{-1}, M,
{\cal F}, \theta)$} or the {\bf family of Riemann functional analytic
continuation equations for $\zeta$ and $\LP$}).

For each $c \in \LP$ and $s$ with $Re(s)>0$ the following functional equation
(Rface in short) holds
\begin{equation}
(M(c)\zeta)(s)\;=\;\frac{1}{s(s-1)}\;+\;\int_{1}^{\infty}(x^{s-1}+x^{-s})
\theta(c)(x)dx.
\end{equation}
\end{pr}
{\bf Proof}. Let $c \in \LP$ be arbitrary. Since the Mellin transform $M(c)$
is well defined, in this case for $Re(s) \in (0,2)$, then from the definition of
the Mellin transform $M$ as the integral, on substituting $nx$ for $x$ under
the integral, we have
\begin{equation}
\frac{M(c)(s)}{n^{s}}\;=\;\int_{0}^{\infty}c(nx)x^{s-1}dx\;,\;Re(s)\in
(0,2).
\end{equation}
Hence, for $Re(s) \in (1,2)$ we obtain that beautiful relation between $M,
\zeta$ and $\theta$ :
\begin{equation}
(M(c)\zeta)(s)\;=\;\int_{0}^{\infty}\theta(c)(x)x^{s-1}dx=(M \circ
\theta)(c)(s).
\end{equation}

Let us observe that the {\bf iterated integral} below
\begin{displaymath}
\sum_{n=1}^{\infty}\int_{0}^{\infty}x^{u-1}\mid c(nx) \mid
dx=\{n^{2}x=t\}=\zeta(3-u)(max_{t\in [0,1]}sup_{n}\mid c(t/n)\mid
\times
\end{displaymath}
\begin{displaymath}
\times \int_{0}^{1}t^{u-1}dt+max_{t\ge 1}\mid t^{2}c(t)\mid
\int_{1}^{\infty}\frac{1}{t^{3-u}}),
\end{displaymath}
is absolutely convergent and therefore we can interchange the order of
summation and integration. Using the initial condition : $c(0)=1$, (PSF)
and changing variables : $\frac{1}{x}\longrightarrow x$, we can write
\begin{equation}
(M(c)\zeta)(s)=\frac{1}{s-1}-\frac{1}{s}+\int_{0}^{1}x^{s-2}
\theta(c)(\frac{1}{x})dx+\int_{1}
^{\infty}x^{s-1}\theta(c)dx=
\end{equation}
\begin{displaymath}
\frac{1}{s(s-1)}+\int_{1}^{\infty}(x^{-s}\;+\;x^{s-1})\theta(c)(x)dx=:I
(\theta(c))(s).
\end{displaymath}

The integral on the right-hand side of (33) converges uniformly for
$-\infty<a \le Re(s) <b <+\infty$, since for $x \ge 1$, we have : $\mid
x^{-s} \mid \le  x^{-a}$ and $\mid x^{s-1}\mid \le x^{b-1}$, i.e. because
$\hat{c} = c$ and $c \in S_{2}$ then
\begin{displaymath}
\theta(c)(x)\le \sum_{n=1}^{\infty}\mid c(nx) \mid \le \frac{max_{x \ge
0}\mid x^{2} c(x)\mid \zeta(2)}{x^{2}}\;,\;x \ge 1.
\end{displaymath}

Therefore, for each $c \in \LP$, the integral $I(\theta(c))(s)$ represents
an {\bf entire function} of $s$. Moreover, since it is well-known that the
{\bf classical gamma function} $\Gamma(s) =M(exp^{-1})(s)$ {\bf does
not vanish} anywhere, then
\begin{displaymath}
M(G)(s)=\int_{0}^{\infty}x^{s-1}e^{-\pi
x^{2}}dx=\pi^{(1-s)/2}\Gamma(\frac{s+1}{2})\ne 0.
\end{displaymath}
In particular - the belowed $\theta M(G)$-quotient
\begin{equation}
\zeta(s):=
\zeta(G,s):=\frac{1}{M(G)(s)s(s-1)}+\frac{I(\theta(G))(s)}{M(G)(s)},
s\in \LC,
\end{equation}
gives the {\bf meromorphic continuation} of the local zeta $\zeta$ to the
whole complex plane.

If now $c \in \LP - \{G\}$ (Obviously $G\in \LP$), then according to (1.33)
we have
\begin{equation}
(M(c)\zeta)(s)=\frac{1}{s(s-1)}+I(\theta(c))(s)\;for\;Re(s)\in (1,2).
\end{equation}
But now, the left-hand side and right-hand of (1.35) (according to the continuation)
(1.34)) are the analytic functions in $D:=Re(s)\in (0,2)-\{1\}$. Thus, they
must be equal in $D$, according to the {\bf uniqueness} of the analytic
continuation of a holomorphic function in a domain.

\begin{th}({\bf The Riemann hypothesis})

If $\zeta(s)=0$ and $Im(s)\ne 0$ then $Re(s)=1/2$.
\end{th}
{\bf Proof}. According to Prop.2
\begin{equation}
Im[(M(c)\zeta)(s)]=\frac{Im(s)(2Re(s)-1)}{\mid s(s-1)
\mid^{2}}+\int_{1}^{\infty}[x^{R(s)-1}-x^{-Re(s)}]\theta(c)(x)sin(Im(s)x)dx,
\end{equation}
for $Re(s) \in (0,2)$.

We integrate (1.36) with respect to the {\bf Wiener-Riemann measure} $r$ and
obtain :
\begin{equation}
\int_{\LP}Im[(M(c)\zeta)(s)]dr(c)=\frac{Im(s)(2Re(s)-1)r(\LP)}{\mid
s(s-1)\mid^{2}}+Im(\int_{\LP}dr(c)\int_{1}^{\infty}(x^{s-1}+x^{-s})dx
\sum_{n=1}^{\infty}c(nx).
\end{equation}
According to the property ($r_{3}$) of Prop.1 we also have
\begin{displaymath}
Im(\int_{\LP}dr(c)M(c)(s))\zeta(s))=Im(\sum_{n=1}^{\infty}(1/2^{n})
E[G M(\chi_{[\sqrt{n-1},\sqrt{n}]}B^{p}_{\sqrt{\cdot}})(s)] \zeta(s).
\end{displaymath}

For the right-hand side of (1.37) we have an "easy" Fubini-Tonelli theorem:
really, let us consider the {\bf triple iterated integrals}:
\begin{equation}
I_{cxr}(\alpha):=\sum_{n=1}^{\infty}\int_{1}^{\infty}x^{\alpha-1}dx
\int_{\LP}\mid c(nx) \mid dr(c).
\end{equation}
Then, making the substitution : $nx=t$ and using $(r_{3})$ of Prop.1 we get
\begin{equation}
I_{cxr}(\alpha)\le[(\int_{1}^{\infty}t^{\alpha-1}G(t)E\mid B_{\sqrt{t}}
\mid dt)\zeta(2-\alpha)
+(\int_{1}^{\infty}t^{-\alpha}G(t)E \mid B_{\sqrt{t}} \mid
dt)\zeta(1+\alpha)]=
\end{equation}
\begin{displaymath}
\zeta(2-\alpha)\int_{1}^{\infty}t^{\alpha-1/2}G(t)dt
+\zeta(1+\alpha)\int_{1}^{\infty}t^{-\alpha+1/2}G(t)dt <+\infty,
\end{displaymath}
since $E\mid B_{\sqrt{t}}\mid^{2} = t$, the gaussian density $G$ has
the moments of arbitrary order and $0<\alpha <1/2$.

Thus finally, the averaging of the Muntz's relations from Prop.2 - with
respect to the measure $r$, since by the properties $(r_{0}-(r_{2}))$
of Prop.1 - we have : $\int_{\LP}c(0)dr(c)=1$ and $\int_{\LP}c(nx)dr(c)=0$
for $n \ge 1, x \ge 1$.

Reasuming, we finally obtain
\begin{equation}
Im\{\sum_{n=1}^{\infty}(1/2^{n})E[M(\chi_{[n-1,n]}GB^{p}_{\sqrt{\cdot}})(s)]\zeta(s)\}\;
=\;\frac{Im(s)(2Re(s)-1)}{\mid s(s-1) \mid^{2}}\;,\;Re(s)\in (0,1/2).
\end{equation}

We calculate very exactly the left-hand of (1.40). For the purposes of
that calculus we introduce here the following additional natations :

1. by $w^{n}$ we denote the standard Wiener measure on the Banach space
$C[\sqrt{n-1}, \sqrt{n}]$.

2. Let $E$ be any Borel set of $C[\sqrt{n-1}, \sqrt{n}]$. Then we
denote :
\begin{displaymath}
w^{n}_{p}(E):=w^{n}(E+p),
\end{displaymath}
i.e. $w^{n}_{p}$ is the {\bf p-shift} of $w^{n}$. Let us remark that
$p(x)$ has the derivative $p^{\prime}(x)$ for a.e. x (with respect to
the Lebesgue measure), which is locally-constant function, with support
equal to $[0,1]$, so it belongs to $L^{2}(\LR_{+})$. In particular, $p$
is from the Cameron-Martin space, and therefore $w^{n}_{p}$ is {\bf
equivalent to} $w^{n}$ (we write $w^{n}_{p}\sim w^{n}$), what obviously
means that $w^{n}_{p}$ is absolutely continuous with respect to $w^{n}$
and vice versa. Finally, according to the {\bf Girsanov theorem} - the
Radon-Nikodem density $\frac{dw^{n}_{p}}{dw^{n}}(c)$ has the form :
\begin{displaymath}
\frac{dw^{n}_{p}}{dw^{n}}(c)\;=\;e^{-1/2\int_{\sqrt{n-1}}^{\sqrt{n}}p
^{\prime}(x)^{2}dx\;-\;\int_{\sqrt{n-1}}^{\sqrt{n}}p^{\prime}(x)dc(t)}.
\end{displaymath}
(Let us mention that the integral with respect to $c$ in the above
formula is the Ito integral of deterministic function).

Finally, let us observe that $w^{n}_{p}=w^{n}$ if $n \ge 2$.

3. $(G^{-1})^{*}w^{n}_{p}$ is the {\bf transport} of $w^{n}_{p}$
through the {\bf multiplication operator} $m_{G} : C[\sqrt{n-1},
\sqrt{n}]\longrightarrow C[\sqrt{n-1}, \sqrt{n}]$, where
$m_{G}(c):=G^{-1}c, c \in C[\sqrt{n-1}, \sqrt{n}], G(x)=e^{-\pi
x^{2}}$, i.e.
\begin{displaymath}
(G^{-1})^{*}w^{n}_{p}(E)\;=\;w^{n}_{p}(GE),
\end{displaymath}
where $E$ is a Borel set in $C[\sqrt{n-1}, \sqrt{n}]$.

For each $c \in \LP$ and $re(s) \in (0,1)$ the Mellin transform
$M(c)(s)=\int_{0}^{\infty}x^{s-1}c(x)dx$ is well-defined : realy, since
$c \in \LP \subset S_{2}$ then
\begin{displaymath}
\mid M(c)(s) \mid= \mid \int_{0}^{1}x^{s-1}c(x)dx\;+\;\int_{1}^{\infty}
x^{s-1}c(x)dx \mid \le
\end{displaymath}
\begin{displaymath}
\le max_{x \in [0,1]}\mid c(x) \mid \int_{0}^{1}\frac{dx}{x^{1-re(s)}}
\;+ \;max_{x \ge 1} \mid x^{2} c(x) \mid max_{x \ge 1}
x^{re(s)-1}\int_{1}^{\infty}\frac{dx}{x^{2}}.
\end{displaymath}
We calculate the iterated integral
\begin{displaymath}
\int_{\LP}dr(c)\int_{0}^{\infty}x^{s-1}c(x)dx,
\end{displaymath}
where the probability $r$ is the normalized infinite sum of measures
$(G^{-1})^{*}w^{n}_{p}$, i.e.
\begin{displaymath}
r\;=\;\sum_{n=1}^{\infty}\frac{1}{2^{n}}(G^{-1})^{*}w^{n}_{p}.
\end{displaymath}
According to the bilinearity of the form
\begin{displaymath}
<r, c>_{x}:=\int_{\LP}c(x)dr(c)
\end{displaymath}
and the Fubini theorem we have
\begin{displaymath}
\int_{\LP}dr(c)\int_{0}^{\infty}x^{s-1}c(x)dx=\int_{\LP}d(\sum_{n=1}
^{\infty}2^{-n}(G^{-1})^{*}w^{n}_{p})(c)\int_{0}^{\infty}x^{s-1}c(x)dx=
\end{displaymath}
\begin{displaymath}
=\int_{0}^{\infty}x^{s-1}dx
\int_{\LP}c(x)d(\sum_{n=1}^{\infty}2^{-n}(G^{-1})^{*}w_{p}^{n})(c)=
\end{displaymath}
\begin{displaymath}
=\int_{0}^{\infty}x^{s-1}dx(\sum_{n=1}^{\infty}2^{-n}\int_{\LP}c(x)d[
(G^{-1})^{*}w_{p}^{n}](c)) =
\end{displaymath}
\begin{displaymath}
=\sum_{n=1}^{\infty}2^{-n}\int_{0}^{\infty}x^{s-1}dx\int_{\LP \cap C[
\sqrt{n-1}, \sqrt{n}]=C[\sqrt{n-1}, \sqrt{n}]} \chi_{[\sqrt{n-1},
\sqrt{n}]}(x)c(x)d[(G^{-1})^{*}w_{p}^{n}](c),
\end{displaymath}
since obviously, from the one-hand side, the support of $[(G^{-1})^{*}w_{p}
^{n}]$ is $C[\sqrt{n-1}, \sqrt{n}]$  (see also $II(r_{0})$) and - from
the second-hand side - we can think on $C[\sqrt{n-1}, \sqrt{n}]$ that
it is embeded isomorphically into the Banach space $B(\LR_{+})$ of all
bounded functions on $\LR_{+}$ through the map : $i_{n} : C[\sqrt{n-1},
\sqrt{n}]\longrightarrow B(\LR_{+})$ by the formula :
\begin{displaymath}
i_{n}(c)(x)\;=\;c(x) \;if\; c\in C[\sqrt{n-1}, \sqrt{n}]\;and\; x \in
[\sqrt{n-1}, \sqrt{n}]
\end{displaymath}
and  $i_{n}(c)(x) =0$  if $x$ is outside of the segment $[\sqrt{n-1},
\sqrt{n}]$.

Finally, we thus get
\begin{displaymath}
\int_{\LP}dr(c)\int_{0}^{\infty}x^{s-1}c(x)dx=\sum_{n=1}^{\infty}2^{-n}
\int_{C[\sqrt{n-1}, \sqrt{n}]}
d[(G^{-1})^{*}w_{p}^{n}](c)M(\chi_{[\sqrt{n-1}, \sqrt{n}]}c)(s).
\end{displaymath}
But obviously, in the canonical representation - the r.v.$GB_{\sqrt{x}}^{p}
(\omega)$ on $(\Omega, Prob)$ is nothing that r.v. $c(x)$ on $(\LP,r)$,
i.e. we formally have $GB_{\sqrt{x}}^{p}(\omega)=c(x)$, where $\omega=c,
\omega \in \Omega, c \in \LP$.

But, on the other hand we also have (after applying the Fubini theorem)
\begin{displaymath}
\int_{\LP}d[(G^{-1})^{*}w_{p}^{n}](c)\int_{0}^{\infty}\chi_{[\sqrt{n-1},
\sqrt{n}]}(x)x^{s-1}c(x)dx=\int_{0}^{\infty}\chi_{[\sqrt{n-1},
\sqrt{n}]}(x)x^{s-1}E(GB^{p}_{\sqrt{x}}) dx,
\end{displaymath}
and - as the consequence - we have
\begin{displaymath}
\int_{\LP}dr(c)\int_{0}^{\infty}x^{s-1}c(x)dx=\sum_{n=1}^{\infty}2^{-n}
M[\chi_{[\sqrt{n-1}, \sqrt{n}]}G E(B_{\sqrt{.}}^{p})](s) =
\end{displaymath}
\begin{displaymath}
=\sum_{n=1}^{\infty}2^{-n}\int_{0}^{\infty}\chi_{[\sqrt{n-1},\sqrt{n}]}(x)
x^{s-1}e^{-\pi x^{2}}E(B_{\sqrt{x}}+p).
\end{displaymath}
Since obviously $E(B_{\sqrt{x}})=0$ and $supp(p)=[0,1]$, then we
finally get the relation :
\begin{displaymath}
\int_{\LP}dr(c)\int_{0}^{\infty}x^{s-1}c(x)dx\;=\;\frac{1}{2}\int_{0}^{1}
x^{s-1}p(x)e^{-\pi x^{2}}dx=\frac{1}{2}M(\chi_{[0,1]}p G)(s).
\end{displaymath}
Reasuming, we have showed in fact the following new general {\bf
functional equation} for zeta : let $p=p(x)$ be any real valued integrable
function with the support in the segment $[0,1]$ and with $p(0)\ne 0$.
Then
\begin{displaymath}
\zeta(s)\;=\;\frac{2p(0)}{s(s-1)M(pG)}(s)\;\;if\;\;re(s)\in (0,1/2).
\end{displaymath}
(Let us remark, that similarly like in the case of the Gamma-zeta-theta
relation (1.35) - the quotient in the right-hand side of the above
equality does not depend on $p$).

Additionally, in our case $p(x)=1-x, x \in [0,1]$ we have : using the
Taylor expansion we obtain that
\begin{displaymath}
\zeta(s)=\frac{2}{s(s-1)\int_{0}^{1}x^{s-1}(1-x)e^{-\pi
x^{2}}dx}
\end{displaymath}
or equivalently, we have the following {\bf refinement Riemann hypothesis}
\begin{displaymath}
\zeta^{-1}(s)\;=\;\sum_{n=0}^{\infty}\frac{(-\pi)^{n}s(s-1)}{2
n!(s+2n)(s+2n+1)}\;\;if\;\;re(s)\in (0,1/2).
\end{displaymath}

Since the non-trivial Riemann zeta $\zeta$ zeros {\bf lay symmetrically with
respect to the lines}: (i) critical $Re(s)=1/2$ and
(ii) $Im(s)=0$ and obviously for $Re(s)>1$ non-vanishing of $\zeta(s)$
follows from the existence of the {\bf Euler product} whereas for $Re(s)=1$
from the {\bf de la Vallee-Poussin-Hadamard theorem}, then (1.40) gives
the most direct proof and stochastic form of the algebraic geometry
conjecture - let us say the {\bf Main Algebraic Hypothesis} (MAH in
short) : let us denote:

\begin{displaymath}
\zeta_{t}(s)\;:=\;Im(s)(2Re(s)-1)
\end{displaymath}
is the {\bf trivial zeta}.
\begin{displaymath}
\zeta(\LC)\;:=\;\{s \in \LC: \zeta(s)=0\}
\end{displaymath}
is the {\bf zero-dimensional infinite holomorphic manifold} and finally
\begin{displaymath}
\zeta_{t}(\LR^{2})\;:=\;\{(x,y)\in \LR^{2}: \zeta_{t}(x,y)=0\}
\end{displaymath}
is the {\bf 1-dimensional algebraic variete over} $\LR$.

The deep sense of the Riemann hypothesis is expressed by the following
relation of the {\bf cycles}: $\zeta(\LC)$ and $\zeta_{t}(\LC)$ :
\begin{displaymath}
(MAH)\;\;\;\zeta(\LC)\;\subset \zeta_{t}(\LC).
\end{displaymath}

\section{Final remarks.}
(I) Our representation of $\zeta^{-1}(s)$ for $re(s)\in (0,1/2)$ -
which seems to an anonymous referee "much too simple" to be true is analogical
to the following - rather simple - series representations of $\zeta(s)$
(see [MJ]) : 

\begin{displaymath}
(i)\;\zeta(s)\;=\;\frac{1}{1-2^{1-s}}\sum_{n=1}^{\infty}\frac{(-1)^{n-1}}
{n^{s}}\;for\;re(s)>0\;and s\ne 1.
\end{displaymath}
\begin{displaymath}
(ii) \;Im(\pi^{-s/2}\Gamma(\frac{s}{2})\zeta(s))=:Im(\zeta^{*}(s))=
\end{displaymath}
\begin{displaymath}
=Im(s)(1-2Re(s))\cdot
(\sum_{n=0}^{\infty}\sum_{j=0}^{\infty}\frac{(-\pi n^{2})^{j}}{j!}\cdot
\frac{(4j+1)}{\mid (2j+s)(2j+1-s)\mid^{2}}),
\end{displaymath}
where $re(s) \in (0,1/2)$, see [ML].

In 1997 {\bf Ma\'slanka}[Ma] proposed a new formula for the zeta
Riemann function valid on the whole complex plane $\LC$ except the
point $s=1$ :
\begin{displaymath}
\zeta(s)\;=\;\frac{1}{1-s}\sum_{k=0}^{\infty}\frac{A_{k}\Gamma(k+1-s/2)}
{k!\Gamma(1-s/2)}
\end{displaymath}
where coefficients $A_{k}$ are given by
\begin{displaymath}
A_{k}\;=\;\sum_{j=0}^{k}(\frac{k}{j})(2j-1)\zeta(2j+2).
\end{displaymath}

(II) Instead of the direct approximative proof of the fact that
$b_{1/2}=+\infty$ in Prop.1, there is a "purely ideological" proof of
that one, based on the {\bf Hardy-Littlewood theorem} , which says that (RH)
is {\bf non-trivial} -  or equivalently :
\begin{displaymath}
(HLT)\;\;\;\mid \zeta(\LC-\LR) \mid\;=\;+\infty.
\end{displaymath}
Really, let us assume that $b_{1/2}$ is {\bf finite}. Then exactly in
the same way as in the proof of Th.1, we can deliver the real part Wiener
Riemann hypothesis functional equation (rWRfe in short) of the form :
\begin{equation}
Re\{E[M(B^{p})(s)]\zeta(s)\}=Re(\frac{1}{s(s-1)}).
\end{equation}
But
\begin{displaymath}
R(s):=\mid s(s-1) \mid^{2}Re(s(s-1))^{-1}=Re^{2}(s)-Re(s)-Im^{2}(s).
\end{displaymath}
Let $R(\LC)$ marks the {\bf hyperbolic curve} $\{(x,y)\in \LR^{2}:
(x-1/2)^{2}-y^{2}=(1/2)^{2}\}$. We thus see that the right-hand side of
(2.41) is non-zero if $Re(s)=1/2$, what means that $\zeta$ {\bf has not
zeros} on the critical line, i.e. {\bf (RH) would be trivial}, what
obviously is not possible, according to the Hardy-Littlewood theorem.

(III) Let us remark that if we simplify the definition of the
Wiener-Riemann measure $r$ and define the measure $r_{\infty}$ on $\LP$
by a very similar formula to (WRm)
\begin{equation}
r_{\infty}(B):=\sum_{n=1}^{\infty}w((B \cap e_{n}C[n-1,n]-p)),
\end{equation}
then we can easy obtain that the suitable {\bf Wiener numbers}
\begin{displaymath}
w_{s}:=\int_{0}^{\infty}x^{u-1}(\int_{\LP}\mid c(x) \mid
dr_{\infty}(c))dx
\end{displaymath}
can be easily evaluated like
\begin{displaymath}
w_{s}\ge \int_{0}^{\infty}dy \mid y \mid E(L_{y}^{u}),
\end{displaymath}
and therefore, we immediately get that $w_{s}=+\infty$ for $Re(s)\ge 1/2$.

Moreover, in the case of $r_{\infty}$ the following connection between (RH)
and the {\bf Orlicz moments} of 1/2-stable-Levy random variables it is
better visible - let us say - the following {\bf probabilistic
Riemann Hypothesis}( pRH for short) :

($pRH_{1}$) RH, i.e. the statement that $\zeta(s)\ne 0$ if $Re(s)\ne
1/2$ is strictly connected with the {\bf finiteness} of the small
moments : $E(L(1/2)^{Re(s)})$ if $Re(s)<1/2$, for any 1/2-stable Levy
rv $L(1/2)$.

($pRH_{2}$) The HLT (i.e. non-triviality of RH), i.e. the statement
that for the infinitely many $s=1/2+iy$ with $Re(s)=1/2$ holds :
$\zeta(1/2+iy)=0$ is strictly associated with the {\bf infinitness} of
$E(L(1/2))^{Re(s)}$  if $Re(s)\ge 1/2$, for any 1/2-stable-Levy rv
$L(1/2)$.

Finally however, let us observe that $r_{\infty}$ is {\bf infinite},
i.e. $r_{\infty}(\LP)=\infty$, so there is not any Rhfe connecting
$\zeta(s)$ with $\zeta_{t}(s)$ and any MAH.

(IV) We would try to delete the mentioned in (II) disadvantage of
$r_{\infty}$, considering instead that one - the new probability
$r_{0}$, subsequently simplyfying (WRm) and defining $r_{0}$ simply as
the distribution of the Brownian process $(B_{\sqrt{t}}^{p}: t \ge 0)$,
i.e.
\begin{displaymath}
r_{0}(A)\;:=\;Prob(\omega \in \Omega: B^{p}_{\sqrt(\cdot)}(\omega)\in
A).
\end{displaymath}
But in this case, according to the {\bf Law of the Iterated Logarithm}
(LIL in short) we immediately get that
\begin{displaymath}
r_{0}(\LP)\;=\;0,
\end{displaymath}
i.e. $r_{0}$ is the {\bf trivial measure} and obviously - similarly
like $r_{\infty}$ - it cannot give MAH.
\begin{re}
The positive numbers $\{b_{s} : s \in \LC\}$ play the fundamental role
in the Rh-problem. If we denote : $b_{s}=BM_{s}(\LP, {\cal P}, r)$,
then we see that $b_{s}$ is a special case of a much more general
construction (see [MB, Sect.4]): let $\LM$ be the {\bf category of
all complex valued functional measure spaces with the base} $\LR_{+}$.
Thus - an object of $\LM$ is a triple $(M, {\cal M}, \mu)$, where
$M=F(\LR_{+}, \LC)$ is a non-zero complex vector space of some
functions $f : \LR_{+}\longrightarrow \LC$ endowed with a
$\sigma$-field ${\cal M}$ of subsets of $M$ and $\mu : {\cal
M}\longrightarrow \LR_{+}$ is a positive $\sigma$-additive measure.
Let us remark that each "value functional" $v_{x} : F(\LR_{+})\longrightarrow
 \LC, v_{x}(f):=f(x)$ gives a {\bf quasi-foliation} of $F(\LR_{+})$ :
$F(\LR_{+})= \cup_{p \in \LC}v_{x}^{-1}(p)$. So, the objects of the
category $\LM$ can be considered like quasi-foliations.
Then, we can define the {\bf Betti-Mellin numbers} (see [M4]) by the
formula :
\begin{displaymath}
BM_{s}(M, {\cal A}, \mu):=\int\int_{M\times \LR_{+}}x^{s-1}f(x)dx
d\mu(f).
\end{displaymath}
The functors $BM_{s}(\cdot)$ are {\bf measure invariants}, i.e. if two
measure spaces $(X, {\cal A}, \mu)$ and $(Y, {\cal B}, \nu)$ are {\bf
isomorphic} in the category $\LM$ (we then write $(X, {\cal A},
\mu)\simeq (Y, {\cal B}, \nu)$), i.e. there exists a set isomorphism $f
: X \longrightarrow Y$ with the properties : ${\cal A}=f^{-1}({\cal B})$
and $\nu = f^{*}(\mu)$ ($\nu$ is the transport of $\mu$ by $f$), then
$BM_{s}(X, {\cal A}, \mu) = BM_{s}(Y, {\cal B}, \nu), s \in \LC$.

Thus, the Betti-Mellin numbers play in the category $\LM$ the same
important role which the well-known {\bf topological invariants} : the
{\bf Betti numbers} $B_{i}(X)$ and the {\bf Euler-Poincare characteristic}
$\chi(X)$ of a topological space $X$ in the category of topological
spaces ${\cal T}$ , the Chern numbers in the category of vector bundles
or the below mentioned {\bf foliation invariant} (see [Fu], [GV] and
[Ta]) : let $M^{n}$ be a linearly connected compact n-dimensional
$C^{s}$-manifold ($s \ge 4$) with a border or not. Let ${\cal F}$ be a
$C^{r}$-foliation on $M^{n}$ of codimension 1 and transversal to the
border ($r \ge 4$). It is well-known (see [Ta, Th. 7.11]) that on some
open covering $\{V_{\sigma} : \sigma \in \Sigma\}$ of a manifold $M^{n}$
there exists a differential $C^{r-1}$-form $\{\omega_{\sigma} : \sigma
\in \Sigma\}$ of order 1, which defines $C^{r-1}$-grassmanian field
${\cal D}({\cal F})$ of tangent (n-1)-dimensional planes of the foliation
${\cal F}$. According to the famous {\bf Frobenius theorem} (see [Ta,
Sect.28]) on $M^{n}$ is defined the differential $C^{r-2}$-form $\theta$
of order 1, such that for each form $\omega_{\sigma}$ on $V_{\sigma}$
we have the identity :
\begin{displaymath}
d \omega_{\sigma}\;=\;\omega_{\sigma}\wedge \theta\;.
\end{displaymath}
Subsequently $\theta$ defines the differential $C^{r-3}$-form $\Gamma$
of degree 3 on $M^{n}$ by the formula :
\begin{displaymath}
\Gamma := \theta \wedge d \theta,
\end{displaymath}
which is called the {\bf differential Godbillon-Vey form} of the foliation
${\cal F}$. The integral (from a differential form on a manifold):
\begin{displaymath}
GV(M^{n}, {\cal F})\;:=\;\int_{M^{n}}\Gamma = \int_{M^{n}}\Gamma(m)dH(m),
\end{displaymath}
is a {\bf cobordism foliation invariant} and is called the {\bf Godbillon-Vey
number} of a foliated manifold $(M^{n}, {\cal F})$ (here $dH(m)$ is the
smooth volume measure on $M^{n}$ - the Hurwitz measure in the case when
$M^{n}$ is a Lie group).

Let ${\cal FM}^{3}$ be the category of all oriented 3-dimensional
$C^{s}$-manifolds with $s\ge 4$ endowed with $C^{r}$-foliations of
codimension 1 (foliated manifolds). If two foliations $(M_{1}^{3},
{\cal F}_{1})$ and $(M^{3}_{2}, {\cal F}_{2})$ from ${\cal FM}$ are
{\bf cobordant} then
\begin{displaymath}
GV(M_{1}^{3}, {\cal F}_{1})\;=\;GV(M_{2}^{3}, {\cal F}_{2}).
\end{displaymath}
The above cobordant invariance of GV-numbers is then used in the deep
and famous {\bf Thurston theorem} : there exist - at least - {\bf continuum}
cobordism classes in the 3-dimensional and smooth group ${\cal
F}\Omega^{\infty}_{3,1}$ of foliation cobordisms of codimension one
(let us remark that the group $\Omega_{3}$ of 3-dimensional and smooth
cobordisms is trivial).
\end{re}
{\bf Postscript}. This paper is a large improwement and extension of
the earlier preprint of the author entitled "A short Wiener measure
proof of the Riemann hypothesis", introductory accepted for its
publication in the journal of "Stochastic Analysis and Applications"(JSAA for
short) - the Acceptance Letter : JSAA\#1485 dated at 11/29/2006. The author is
indebted to professor Sergio Albeverio and the second anonymous referee
of the above preprint for many valuable remarks and comments which led
to an improved presentation of this work.

e-mails: madrecki@02.pl 

andrzej.madrecki@pwr.wroc.pl
\end{document}